\documentclass[11pt]{article}
\usepackage{amssymb,amsfonts,amsmath,latexsym,epsf,tikz,url}

\newtheorem{theorem}{Theorem}[section]
\newtheorem{proposition}[theorem]{Proposition}

\newtheorem{corollary}[theorem]{Corollary}
\newtheorem{lemma}[theorem]{Lemma}

\newtheorem{definition}[theorem]{Definition}

\newcommand{\proof}{\noindent{\bf Proof.\ }}
\newcommand{\qed}{\hfill $\square$\medskip}

\begin{document}

\title{Distinguishing number and distinguishing index of lexicographic  product of two graphs}

\author{
Saeid Alikhani  $^{}$\footnote{Corresponding author}
\and
Samaneh Soltani
}

\date{\today}

\maketitle

\begin{center}
Department of Mathematics, Yazd University, 89195-741, Yazd, Iran\\
{\tt alikhani@yazd.ac.ir, s.soltani1979@gmail.com}
\end{center}

\begin{abstract}
	The distinguishing number (index) $D(G)$ ($D'(G)$) of a graph $G$ is the least integer $d$ such that $G$ has an vertex labeling (edge labeling)  with $d$ labels  that is preserved only by a trivial automorphism. The lexicographic product of two graphs  $G$ and $H$, $G[H]$ can be obtained from $G$ by substituting a copy
	$H_u$ 	of $H$ for every vertex $u$ of $G$ and then joining all vertices of $H_u$ with all vertices of $H_v$ 	if $uv\in E(G)$. 
In this paper we obtain some sharp bounds for  the distinguishing number and the distinguishing index 
of lexicographic product of two graphs. As consequences, we prove that if $G$ is a connected graph with a special condition on automorphism group of $G[G]$ and  $D(G)> 1$,  then for every natural $k$, $D(G)\leq D(G^k)\leq D(G)+k-1$, where $G^k=G[G[...]]$. Also we prove that  all lexicographic powers of $G$, $G^k$ ($k\geq 2$) can be distinguished by
at most two edge labels. 
\end{abstract}


\medskip
\noindent{\bf AMS Subj.\ Class.:} 05C15, 05E18

\section{Introduction}
Let $G = (V ,E)$ be a simple  graph with $n$ vertices.   Throughout this paper we consider only simple graphs.
The set of all automorphisms of $G$, with the operation of composition of permutations, is a permutation group
on $V$ and is denoted by $Aut(G)$.  
A labeling of $G$, $\phi : V \rightarrow \{1, 2, \ldots , r\}$, is  $r$-distinguishing, 
if no non-trivial  automorphism of $G$ preserves all of the vertex labels.
In other words,  $\phi$ is $r$-distinguishing if for every non-trivial $\sigma \in Aut(G)$, there
exists $x$ in $V$ such that $\phi(x) \neq \phi(x\sigma)$. 
The distinguishing number of a graph $G$ has defined by Albertson and Collins \cite{Albert} and  is the minimum number $r$ such that $G$ has a labeling that is $r$-distinguishing.  
Similar to this definition, Kalinkowski and Pil\'sniak \cite{Kali1} have defined the distinguishing index $D'(G)$ of $G$ which is  the least integer $d$
such that $G$ has an edge colouring   with $d$ colours that is preserved only by a trivial
automorphism.  These indices  have  developed  and number of papers published on this subject (see, for example \cite{soltani2,chan,Klavzar,fish}).
For every vertex $v\in V$, the open neighborhood of $v$ is the
set $N_G(v) =\{u\in  V : uv \in  E\}$  and the closed neighborhood is the set $N_G[v] = N_G(v)\cup \{v\}$. 

For two graphs $G$ and $H$, let $G[H]$ be the graph with
vertex set $V(G)\times V(H)$,  such that the vertex $(a,x)$ is
adjacent to vertex $(b, y)$ if and only if $a$ is adjacent to $b$
(in $G$) or $a = b$ and $x$ is adjacent to $y$ (in $H$). The graph
$G[H]$ is the lexicographic product  of $G$ and
$H$.
This product was introduced as the composition of graphs by Harary \cite{Harary}.
The lexicographic product is also known as graph
substitution, a name that bears witness to the fact that $G[H]$ can be obtained from $G$ by
substituting a copy $H_u$ of $H$ for every vertex $u$ of $G$ and then joining all vertices of $H_u$ with
all vertices of $H_v$ if $uv \in E(G)$. For example $K_2[K_3] = K_6$.
It can be seen that the number of edges of $G[H]$ is $\vert  V(G)\vert \vert E(H)\vert + \vert E(G)\vert \vert V(H)\vert^2 $. Also the degree of an arbitrary vertex $(g,h)$ of $G[H]$ is $deg_H h+\vert V(H)\vert deg_G g $.
The distinguishing number and the distinguishing index of some operations of two graphs, such as Cartesian product and corona product have been studied in \cite{soltani2,Klavzar}. Kla\v{v}zar and  Zhu in \cite{Klavzar} have shown that the  Cartesian powers of graphs can be distinguished by two labels. 
In this paper we shall study the distinguishing number and the distinguishing index of lexicographic product of two graphs.  To do this,  we consider the automorphisms of $G[H]$ in this section.  
In Section 2, we study the distinguishing number of $G[H]$. In Section 3, we study the distinguishing index of lexicographic product of two graphs. Here we state some properties of automorphisms of $G[H]$.

Let $\beta$ be an automorphism of $H$, and $(g,h)$ a vertex of $G[H]$. The  permutation of $V (G [H])$ that maps $(g, h)$ into $(g, \beta h)$ and is the identity elsewhere, clearly is in $Aut(G[H])$. Also, if $\alpha \in  Aut(G)$, then the mapping $(g, h) \mapsto (\alpha g, h)$ is an automorphism of $G[H]$.
The group generated by such elements is known as the wreath product $Aut(G)[Aut(H)]$.
Evidently all its elements can be written in the form $(g, h) \mapsto (\alpha g, \beta_{\alpha g}h)$,
where $\alpha$ is an automorphism of $G$ and  $\beta_{\alpha g}$ are automorphisms of $H$. As the example of $K_2[K_2]$ shows, $Aut(G)[Aut(H)]$ can be a proper subgroup of $Aut(G[H])$.  In fact the elements of  $Aut(G)[Aut(H)]$ are the automorphisms that they map the copies of $H$ to each other, completely. The next theorem describes when  $Aut(G)[Aut(H)]$  is equal to $Aut(G[H])$. For the statement of the theorem, we use the relations $S$ and $R$ that are defined as follows:

\begin{definition}{\rm \cite{Sabidussi}}\label{relations}
Let $G$ be a graph. The equivalence relation $R$ and $S$ are defined on $V(G)$ as follows:
\begin{equation*}
g_1 R g_2 \Longleftrightarrow N_G(g_1)=N_G(g_2),~~~g_1 S g_2 \Longleftrightarrow N_G[g_1]=N_G[g_2].
\end{equation*}
\end{definition}

 \begin{theorem}{\rm \cite{Sabidussi}}\label{Gert}
Let $G, H$ be two graphs and  $R, S$ be the relations on $V(G)$ in Definition \ref{relations}. Then a necessary and sufficient condition that $Aut(G[H])=Aut(G)[Aut(H)]$ is that $H$ be connected if $R\neq \Delta$, and that $\overline{H}$ (the complement of $H$) be   connected if $S\neq \Delta$, where $\Delta=\{(g,g) \vert~ g\in V(G)\}$.
 \end{theorem}

We now reply to the question: What is the $Aut (G[H])$ when $\overline{H}$ is  disconnected graph and $G$ has nontrivial automorphism? In \cite{Bird},  Bird, et al. 	have replied to this question for partially ordered sets $G$ and $H$, where $Aut(G)$ consists of all permutations on $G$ that preserve order (and have order preserving inverses). By considering $\overline{G[H]}=\overline{G}[\overline{H}]$, it can be replied to this question exactly the same as Bird, et al.  as follows:

\begin{theorem}\label{thmbird}
Let $G$ and $H$ be two graphs, the  $S$-equivalent pairs in $G$ denoted by $[g_{j1},g_{j2}]$ for $j=1,\ldots , \theta$ and  the connected component of $\overline{H}$ be denoted by $\overline{H}_i$. Also consider the following elements:  
\begin{equation*}
S(ij)=\left\{
\begin{array}{ll}
(g,h)\mapsto (g,h) & h\in \overline{H}_i,\\
(g_{j1},h)\mapsto (g_{j2},h) & h\notin \overline{H}_i,\\
(g_{j_2},h)\mapsto (g_{j_1},h) & h\notin \overline{H}_i,\\
(g,h)\mapsto (g,h) & h\notin \overline{H}_i, g\neq g_{j1},g_{j2}.\\
\end{array}\right.
\end{equation*}
Then the automorphism group of lexicographic product of $G$ and $H$,  $Aut(G[H])$ is the group generated by adding the elements $S(ij)$ to the wreath product $Aut(G)[Aut(H)]$. 
\end{theorem}

\section{The distinguishing number of $G[H]$}
In this section we study the distinguishing number of lexicographic product of two graphs $G$ and $H$. The following theorem gives  sharp bounds  for the distinguishing number of $G[H]$.

\begin{theorem}\label{thmdis1}
Let $G$ and $H$ be two connected graphs, then 
\begin{equation*}
D(H)\leqslant D(G[H])\leqslant D(G)\times D(H).
\end{equation*}
\end{theorem}
\proof
First we prove that $D(H)\leqslant D(G[H])$. By contradiction, we suppose that $D(H) > D(G[H])$. So  in the distinguishing labeling of $G[H]$ with $D(G[H])$ labels, it can be seen that all copies of $H$ have been labeled with less than $D(H)$ labels. Hence for each copy of $H$ there exists a nontrivial automorphism $\beta_g$ of $H$  such that $\beta_g$ do not preserve the labeling of that copy of $H$ in the distinguishing labeling of $G[H]$. So there exists  the following nontrivial automorphism $f$ of $G[H]$ 
\begin{equation*}
f: V(G[H])\rightarrow V(G[H])~with~  f(g,h)= (g,\beta_g h),
\end{equation*}
 such that $f$ do not preserve the labeling of $G[H]$, which is a contradiction.

Now we want to show that $D(G[H])\leqslant D(G)\times D(H)$. For this purpose, we label the vertices of $i$th copy of $H$  with the labels $\{1+(i-1)D(H),2+(i-1)D(H),\ldots ,D(H)+(i-1)D(H) \}$ in a distinguishing way, where $1\leqslant i \leqslant |V(G)|$. This labeling is a distinguishing labeling of $G[H]$, because if $f$ is an automorphism of $G[H]$ preserving the labeling, then with respect to the labeling of copies of $H$, the map $f$ maps each copy of $H$ to itself, and since we labeled each copy of $H$ in a distinguishing way, $f$ is the identity automorphism. Since we used $D(G)\times D(H)$  labels for this labeling, the result follows. \qed

The bounds of $D(G[H])$ in Theorem \ref{thmdis1} are sharp. For the upper bound it is sufficient to consider  the complete graphs $K_n$ and $K_m$, as two graphs  $G$ and $H$, respectively.  Because $K_n[K_m]=K_{nm}$. For the lower bound we consider $G=K_1$, then $G[H]=H$, and so $D(H)= D(G[H])$.
\medskip

If $Aut(G[H])=Aut(G)[Aut(H)]$, then we can improve the upper bound of $D(G[H])$ in Theorem \ref{thmdis1} as follows:
\begin{theorem}\label{thmdis2}
Let $G$ and $H$ be two connected graphs with  $Aut(G[H])=Aut(G)[Aut(H)]$. Then  $ D(H)\leqslant D(G[H])\leqslant  D(H)+M$,
 where $M= min\left \{k: \sum_{m=0}^{k}y_m  \geqslant D(G)\right\}$ and  $$ y_m=\left\{
 \begin{array}{ll}
 1&m=0,\\
 D(H)& m=1,\\
 D(H)+\sum_{i=1}^{m-1}{m-1 \choose i}{D(H)\choose i+1}& m\geqslant 2.
 \end{array}\right.$$
\end{theorem}
\proof
The inequality  $D(H)\leqslant D(G[H])$ has proved in Theorem \ref{thmdis1}. For obtaining the upper bound, we partition the vertices of $G$ by a distinguishing labeling of $G$, i.e., we partition the vertices of $G$ into $D(G)$ classes, say $[1],\ldots , [D(G)]$  such that $i$th class contains the vertices of $G$ having the label $i$, in the distinguishing labeling of $G$, where $1\leqslant i \leqslant D(G)$. By this partition we label the copies of $H$ as follows: First we label the vertices of $H$ with $D(H)$ labels in a distinguishing way, next we do the following changes on the labeling of $H$. Before starting the labeling of the copies of $H$, we introduce the notation $H^{[i]}$ for the set of copies of $H$ corresponding to the elements of $i$th class,  where $1\leqslant i \leqslant D(G)$. In fact we partition the copies of $H$ into $D(G)$ classes such that   $H^{[i]}$ is the symbol of $i$th class. Now we present the labeling of $G[H]$ by the following steps:
  
 Step 1) We label all vertices of  the copies of $H$ that are in $H^{[1]}$, exactly the same as the distinguishing labeling of vertices of  $H$.
 
 Step 2) For the labeling of the vertices of copies in $H^{[i]}$, where $2\leqslant i \leqslant D(H)+1$, we use of the new label $D(H)+1$ in such a way that the label $i-1$ in the all elements of $H^{[i]}$ is replaced by the new label $D(H)+1$,   where $2\leqslant i \leqslant D(H)+1$.
 
 Step 3)  For the labeling of the vertices of the copies in $H^{[i]}$, where $D(H)+2\leqslant i \leqslant 2D(H)+1$, we do the same work as Step 2, with the new label $D(H)+2$, instead of the labels $D(H)+1$.
 
 Step 4) By choosing two labels among the labels $\{1,\ldots , D(H)\}$, and replacing them by the two new labels $D(H)+1$ and $D(H)+2$, we can label the elements of ${D(H) \choose 2}$ other classes of the classes  $H^{[i]}$.
 
 Step 5) We do the same work as Step 2 with the new label $D(H)+3$ instead of labels $D(H)+1$. Next we label $2{D(H) \choose 2}$ other classes   $H^{[i]}$, with the two new labels $D(H)+1$ and $D(H)+3$, also with the labels $D(H)+2$ and $D(H)+3$, exactly the same as Step 4.
 
 Step 6) Now we choose three labels among the labels $\{1,\ldots , D(H)\}$, and replace them by the three new labels $D(H)+1$, $D(H)+2$ and $D(H)+3$.
 
 By continuing this method we obtain that the number of classes can be labeled with the labels $1,\ldots , D(H)+m$, $m\geqslant 1$, such that the label $D(H)+m$ is used in the labeling of each element of classes, is $y_m$ where
 
 \begin{equation*}
 y_m=\left\{
 \begin{array}{ll}
 1&m=0,\\
 D(H)& m=1,\\
 D(H)+\sum_{i=1}^{m-1}{m-1 \choose i}{D(H)\choose i+1}& m\geqslant 2.
 \end{array}\right.
 \end{equation*}
  
  Therefore the number of labels that have been used for the labeling of vertices of all copies of $H$, is $D(H)+M$ where $M= min\left \{k: \sum_{m=0}^{k}y_m  \geqslant D(G)\right\}$. 
   This labeling is a distinguishing vertex labeling of $G[H]$, because if $f$ is an automorphism of $G[H]$ preserving the labeling, then  since $Aut(G[H])=Aut(G)[Aut(H)]$, we have $f(g,h)=(\alpha g,\beta_{\alpha g}h)$, for some automorphism $\alpha$ of $G$ and $\beta_{\alpha g}$ of $H$. With respect to the labeling of copies of $H$, it can be concluded that   $\alpha$ is the identity automorphism on $G$. Since each copy of $H$ have been labeled in a distinguishing way, $\beta_{\alpha g}$  is the identity automorphism on $H$, and so $f$ is the identity automorphism on $G[H]$.\qed
 
\begin{figure}[ht]
\begin{center}
\includegraphics[width=0.7\textwidth]{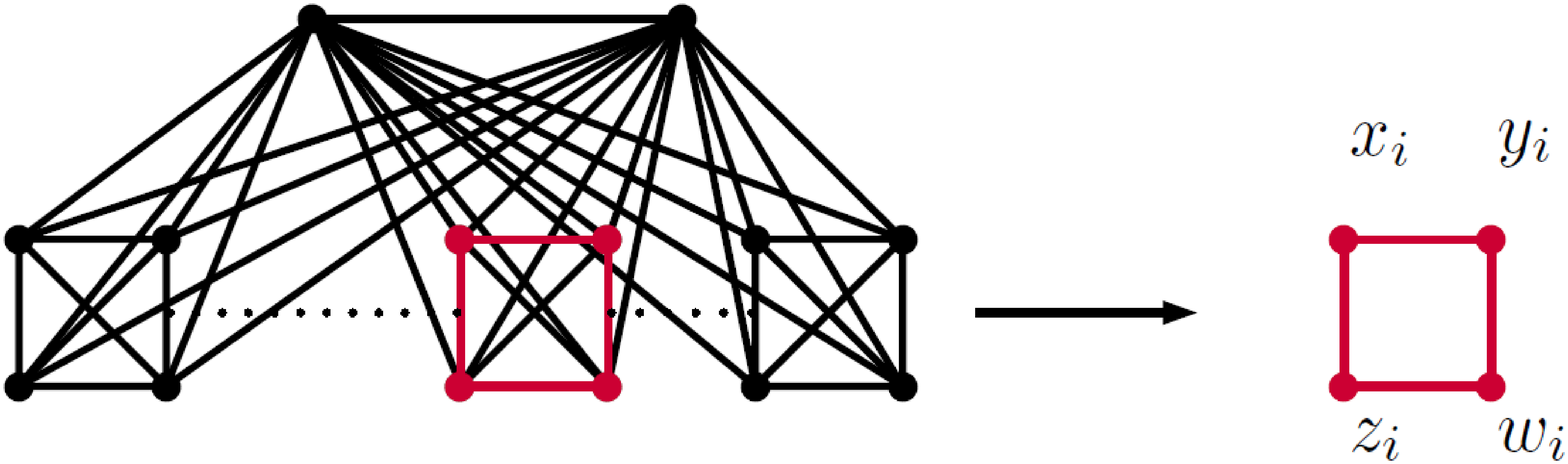}
	\caption{\label{fig1} The place of labels  $x_i,y_i,z_i,w_i$ in $G_n[K_2]$.}
\end{center}
\end{figure}

\begin{figure}[ht]
\begin{center}
\includegraphics[width=1\textwidth]{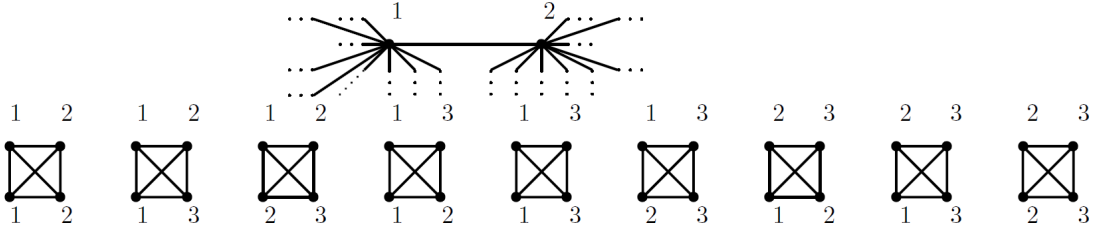}
	\caption{\label{fig2} The $3$-distinguishing labeling of $G_n[K_2]$.}
\end{center}
\end{figure}

 Here we shall show that the upper bound of $D(G[H])$ in Theorem \ref{thmdis2} is sharp. To do this, suppose that   $G_n$ ($n\geq 3$) is a spider graph which has formed by subdividing all of the edges of a star $K_{1,n}$. We state and prove the following lemma: 
 
 \begin{lemma}\label{Gn[K_2]}
For every $n\geq 3$,  $D(G_n[K_2])= \lceil \frac{1+\sqrt{1+8\sqrt{n}}}{2} \rceil$. 
  \end{lemma}
 \proof 
 In an $r$-distinguishing labeling of $G_n$, each of the pairs consisting of a noncentral-nonpendant vertex of a branch of $G_n$ and its pendant neighbor must have different ordered pair of labels. There are $r^2$ possible ordered pairs of labels using $r$ labels, hence $D(G_n)=\lceil \sqrt{n} \rceil$. It is easy to check that $Aut (G_n[K_2])=(Aut (G_n))[Aut (K_2)]$, by Theorem \ref{Gert}. Let $L=\{(x_i,y_i,z_i,w_i)\vert~ 1\leqslant i \leqslant n, x_i,y_i,z_i,w_i \in \mathbb{N}\}$ be  a labeling of the vertices $G_n[K_2]$ except its central vertices (see Figure \ref{fig1}).
  If $L$ is a distinguishing labeling then the label of two central vertices of $G_n[K_2]$ must be different. In addition, the following conditions must satisfy:
 \begin{itemize}
   \item[(i)] $x_i\neq y_i$ and $w_i\neq z_i$, for all $i=1,\ldots , n$.
   \item[(ii)] $\{x_i,y_i,z_i,w_i\}\neq \{x_j,y_j,z_j,w_j\}$, for all $i,j\in \{1,\ldots , n\}$ where  $i\neq j$.
   \end{itemize}
 So there are ${r \choose 2}{r \choose 2}$ possible $4$-arrays of labels using $r$ labels such that they satisfy (i) and (ii), hence $D(G_n[K_2])= \lceil \frac{1+\sqrt{1+8\sqrt{n}}}{2} \rceil$ (see Figure \ref{fig2} for a $3$-distinguishing labeling of $G_n[K_2]$ (note that we do not sketch some edges for blinding clarity)). \qed
 
 Now by Lemma \ref{Gn[K_2]} we see that   for   $n=50$, we have $D(G_{50})=8$, $D(G_{50}[K_2])=5$ and $M=3$. Hence there exists $n$ such that the graph $G_n[K_2]$ obtains  the upper bound of Theorem \ref{thmdis2}.

 \medskip
 As a corollary of  Theorem \ref{thmdis2} we would like to present bounds for the distinguishing number of $G^k=G[G[...]]$.   
   \begin{corollary}
   Let $G$ be a connected graph such that $Aut(G[G])=Aut(G)[Aut G]$, then
   \begin{itemize}
   \item[(i)] If $D(G)> 1$, then $D(G)\leqslant D(G^k)\leqslant D(G)+k-1$.
   \item[(ii)] If $D(G)= 1$, then $D(G^k)=1$.
   \end{itemize}
   \end{corollary}
   \proof
   By Theorem \ref{Gert} it can be seen that if $Aut(G[G])=Aut(G)[Aut (G)]$, then
 $Aut(G^k)=Aut(G)[Aut(G^{k-1})]$ for $k\geq 2$.
       \begin{enumerate}
    \item[(i)] Proof is by induction on $k$. For $k=2$, we observe that by Theorem \ref{thmdis2} the value of $M$ is one, and so the result follows. 
   
   \item[(ii)] if $D(G)= 1$, we observe that by Theorem \ref{thmdis2} the value of $M$ is zero, and so the proof is complete. \qed
   \end{enumerate} 

\section{The distinguishing index of $G[H]$}

In this section we shall study the  distinguishing index of lexicographic  product of two graphs. We begin with the following theorem: 

\begin{theorem}\label{thmindex}
Let $G$ and $H$ be two connected graphs such that $H\neq K_2$ and $Aut(G[H])=Aut(G)[Aut(H)]$. Then $D'(G[H])\leqslant max \{D'(G), D'(H)\}$. 
\end{theorem}
\proof
First we partition the edge set of $G$ into $D'(G)$ classes, say $[1],\ldots, [D'(G)]$, by a distinguishing edge labeling of $G$. In fact, the $i$th class contains the edges of $G$ with the label $i$  in the distinguishing edge labeling of $G$, where  $1\leqslant i\leqslant D'(G)$. For labeling of $G[H]$, we label the edge set of each copy of $H$ in a distinguishing way with $D'(H)$ labels. By the definition of $G[H]$  we know that each edge of $G$, such as $e$ is relaced by the edges which join the corresponding two copies of $H$, in $G[H]$. We denote the set of these replacement edges by $\mathbf{E}$. Now we assign all edges in $\mathbf{E}$, the same label of the edge $e$ in the distinguishing edge labeling of $G$.  This labeling is a distinguishing labeling of $G[H]$, because if $f$ is an automorphism of $G[H]$ preserving the labeling, then  since $Aut(G[H])=Aut(G)[Aut(H)]$, we have $f(g,h)=(\alpha g,\beta_{\alpha g}h)$, for some automorphism $\alpha$ of $G$ and $\beta_{\alpha g}$ of $H$. With respect to the labeling of edges in the set $\mathbf{E}$, where $e\in E(G)$, it can be concluded that   $\alpha$ is the identity automorphism on $G$. Since each copy of $H$ have been labeled in a distinguishing way, $\beta_{\alpha g}$  is the identity automorphism on $H$, and so $f$ is the identity automorphism on $G[H]$.\qed

By Theorem \ref{thmindex}, the lexicographic  product $P_m[P_n]$ of two path of orders $m$ and $n> 2$, respectively, has the distinguishing index equal to $2$, unless $m=2$ (since $Aut(P_2[P_n])\neq Aut(P_2)[Aut(P_n)]$). For the lexicographic product of a cycle $C_n$ with a path $P_m$ we also have $D'(P_m[C_n])=2$ where $m>2$ and $n>5$. The lexicographic  product of two cycles $C_n$ and $C_m$ also has the distinguishing index equal to two, where $n,m > 5$. It is worth noting that these results do not depend on the relation between $n$ and $m$.

\begin{proposition}
If $H$ is a connected graph, then $D'(K_2[H])=\left\{ \begin{array}{ll}
1& \vert V(H)\vert =1,\\
3& \vert V(H)\vert =2,\\
2& \vert V(H)\vert \geqslant 3.
\end{array}\right.$
\end{proposition}
\proof
If $\vert V(H)\vert =1$ then $K_2[H]=K_2$, and so $D'(K_2[H])=1$. If $\vert V(H)\vert =2$ then $K_2[H]=K_4$, and so $D'(K_2[H])=3$. Let $\vert V(H)\vert \geqslant 3$. Since the graph $K_2[H]$ has a nontrivial automorphism, $D'(K_2[H])\geqslant 2$. Now we present a distinguishing edge labeling of $K_2[H]$ with two labels. First we label all edges of the first copy of $H$ with label $1$, and all edges of the second copy of $H$ with the label $2$. Let $V(K_2)=\{x_1,x_2\}$ and  $V(H)=\{y_1,\ldots , y_n\}$ where $n\geqslant 3$. We label the edges $(x_1,y_j)(x_2,y_1),\ldots, (x_1,y_j)(x_2,y_n)$ with $j-1$ labels $2$ and $n-(j-1)$ labels $1$, where $1\leqslant j \leqslant n$. By Theorem \ref{thmbird} this labeling is distinguishing, and so   $D'(K_2[H])=2$.\qed

In \cite{Gorzkowska}, Gorzkowska, et al.
  have obtained the distinguishing index of Cartesian product $K_{1,n}\square P_m$ and $K_{1,n}\square C_m$. We use their method of proof to obtain an upper bound for  distinguishing index of lexicographic 
 product $K_{1,n}[H]$ where $H$ is a graph of order $m\geqslant 2$.
\begin{proposition}
If $H$ is a connected graph of order  $m\geqslant 2$ and $K_{1,n}$ is the star graph with $n\geqslant 2$, then $2 \leqslant D'(K_{1,n}[H])\leqslant max\{ D'(H), \lceil \sqrt[m^2]{n} \rceil\}$, unless $m=2$ and $n=r^4$ for some integer $r$. In the latter case, $2 \leqslant D'(K_{1,n}[P_2])\leqslant  \sqrt[4]{n}+1$.
\end{proposition}
\proof
Since the graph $K_{1,n}[H]$ has a nontrivial automorphism, so $D'(K_{1,n}[H])\geqslant 2$. Now we present a distinguishing edge labeling of $K_{1,n}[H]$. First we label the  edges of each copy of $H$ with $D'(H)$ labels in a distinguishing way.  Let $d$ be a positive integer such that $(d-1)^{m^2}< n \leqslant d^{m^2}$. Denote by $x_0$ the central vertex of the star $K_{1,n}$, by $x_1,\ldots ,x_n$ its pendant vertices, and by   $y_1,\ldots , y_m$ vertices of $H$ where $m \geqslant 2$. Suppose first that $m\geqslant 3$. By Theorem \ref{Gert} every automorphism of $K_{1,n}[H]$ is of the form $f(x,y)=(\alpha x,\beta_{\alpha x}y)$ where $\alpha$ is an automorphism of $K_{1,n}$ and $\beta_{\alpha x}$ an automorphism of $H$. Since we labeled the edges of each copy of $H$ in a distinguishing way, $\beta_{\alpha x}$ is the identity automorphism of $H$ where $f$ is the automorphism of $K_{1,n}[H]$ preserving the labeling. 

We want to show that the remaining edges of $K_{1,n}[H]$ can be labeleded  such that the copies of $H$ also cannot be interchanged. Then the identity automorphism is the only automorphism of $K_{1,n}[H]$ preserving the labeling. A labeling of all edges  yet unlabeled  can be fully described by defining a matrix $L$ with $m^2$ rows and $n$ columns such that in the $j$th column the initial $m$ elements are labels of the  edges
$(x_0,y_1)(x_j,y_1),\ldots ,(x_0,y_1)(x_j,y_m)$,
and the next  $m$ elements are labels of the edges $(x_0,y_2)(x_j,y_1),\ldots ,(x_0,y_2)(x_j,y_m)$, and finally, the last $m$ elements are labels of the edges 
 $(x_0,y_m)(x_j,y_1),\ldots ,(x_0,y_m)(x_j,y_m)$.  If matrix  $L$ contains at least two identical columns, then there exists a permutation of copies of $H$   preserving the labeling, and vice versa.  There are exactly $d^{m^2}$ sequences of length $m^2$ with elements from the set $\{1,\ldots, d\}$, hence
there exists a labeling with $d$ colours such that every column of $L$ is distinct. Therefore,
$D'(K_{1,n}[H])\leqslant max\{D'(H),d\} =max\{D'(H),\lceil \sqrt[m^2]{n} \rceil\}$.

For $m = 2$, we label the edges of $K_{1,n}[P_2]$ in the same way. The only difference is
that each copy of $P_2$  has only one edge, hence the two copies of $P_2$ need not be fixed. This is the case when $n = d^4$, because then each element of $\{1,\ldots, d\}^4$ is a column in $L$, and there exists a permutation of columns of $L$ which together with the transposition of rows
of $L$ defines a non-trivial automorphism of $K_{1,n}[P_2]$ preserving the colouring. Thus we
need an additional label for one edge in a copy of $P_2$. When $n < d^4$, we put the sequence
$(1,1, 1, 2)$ as the first column of $L$, and we do not use the sequence $(1,1, 2, 1)$ any more, thus this labeling breaks the transposition of the rows of $L$, and so all automorphisms of $K_{1,n}[P_2]$.\qed

The following proposition implies that the lexicographic product of $P_n$ ($n\geq 3$) with any connected graph, can be distinguished by two edge labels. 
\begin{proposition}
Let $P_n$ be the path of order $n\geqslant 3$ and $H$ be a connected graph of order $m\geqslant 1$. Then $D'(P_n[H])=2$.
\end{proposition}

\proof
Since the graph $P_n[H]$ has a nontrivial automorphism, so $D'(P_n[H])\geqslant 2$. If $m=1$, then $P_n[H]=P_n$, and so $D'(P_n[H])=2$. Let $m\geqslant 2$, we present a $2$-distinguishing labeling for $P_n[H]$ as follows: We label all edges of each copy of $H$ with the label $1$. If we denote the consecutive vertices of $P_n$ by $x_1,\ldots , x_n$ and vertices of $H$ by $h_1,\ldots h_m$, then for every $1\leqslant i \leqslant n-2$, we label the edges $(x_i,h_j)(x_{i+1},h_1),\ldots ,(x_i,h_j)(x_{i+1},h_m)$ with $j-1$ labels $2$ and $m-(j-1)$ labels $1$ for $1 \leqslant j \leqslant m$. We label the edges $(x_{n-1},h_j)(x_{n},h_1),\ldots ,(x_{n-1},h_j)(x_{n},h_m)$ with $j-1$ labels $1$ and $m-(j-1)$ labels $2$ for $1 \leqslant j \leqslant m$. By Theorem \ref{Gert}, $Aut(P_n[H])=Aut(P_n)[Aut(H)]$, and so the labeling is distinguishing, because if $f$ is an automorphism of $P_n[H]$ preserving the labeling, then $f(x,h)=(\alpha x,\beta_{\alpha x}h)$, for some automorphism $\alpha$ of $P_n$ and $\beta_{\alpha x}$ of $H$. With respect to the labeling of edges between copies of $H$ it is concluded that   $\beta_{\alpha x}$ is the identity automorphism on $H$. Regarding to the labeling of the edges between the first and the second copies of $H$ and the labeling of the edges between the $(n-1)$-th and the last copies of $H$, it follows that $\alpha$   is the identity automorphism of $P_n$. Therefore $f$ is the identity automorphism of  $P_n[H]$.\qed

The following theorem gives an upper bound for the distinguishing index of $G[P_2]$:
\begin{theorem}
Let $G$ be a connected graph such that $Aut(G[P_2])=(Aut (G))[Aut(P_2)]$, then 
\begin{equation*}
D'(G[P_2])\leqslant min \Big\{k: \sum_{m=2}^{k}\left( 2{m-1 \choose 1} +m{m-1 \choose 2}+{m-1\choose 3} \right)\geqslant D'(G)\Big\}.
\end{equation*}
\end{theorem}
\proof
First we partition the edge set of $G$ by a distinguishing labeling into $D'(G)$ classes, say $[1],\ldots , [D'(G)]$ such that $i$th class contains the edges of $G$ having label $i$ in the distinguishing labeling of $G$. Let $[i]=\{e_{i1},\ldots , e_{is_i}\}$ such that $s_i$ is the size of $i$th class, where $1 \leqslant i \leqslant D'(G)$. So each of $e_{i1},\ldots , e_{is_i}$ is replaced by four edges in $G[P_2]$. We denote the set of four edges corresponding to the edge $e_{ij}$ of $G$ by the symbol $\mathbf{E}_{ij}$. For labeling the edges of $G[P_2]$ we first label all copies of $P_2$ with the label $1$. We continue the labeling by the following steps:

\medskip
Step 1) For every $1\leqslant j \leqslant s_1$, we label the edges in $\mathbf{E} _{1j}$  with three labels $1$ and one label  $2$.

\medskip
Step 2) For every $1\leqslant j \leqslant s_2$, we label the edges in $\mathbf{E}_{2j}$  with three labels $2$ and one label  $1$.

So we labeled the corresponding edges to the edges in the first and second classes of $G$ with labels $1$ and $2$.

\medskip
Step 3) For every $1\leqslant j \leqslant s_3$, we do the same work as Step $1$ for labeling the edges in  $\mathbf{E}_{3j}$  with the labels $1$ and  $3$. Also for every $1\leqslant j \leqslant s_4$, we do the same work as Step $2$ for labeling the edges in  $\mathbf{E}_{4j}$  with the labels $1$ and  $3$.

\medskip
Step 4) For every $1\leqslant j \leqslant s_5$ and $1\leqslant j \leqslant s_6$, 
we do the same work as Steps  $1$ and $2$, respectively, with the labels $2$ and  $3$.

\medskip
Step 5)  For every $1\leqslant j \leqslant s_7$, we label four edges in $\mathbf{E}_{7j}$  with the labels $1,2,3,1$.  For every $1\leqslant j \leqslant s_8$, we label four edges in $\mathbf{E}_{8j}$  with the labels $1,2,3,2$.  For every $1\leqslant j \leqslant s_9$, we label four edges in $\mathbf{E}_{9j}$  with the labels $1,2,3,3$.

So we labeled the corresponding edges to the classes $[3],\ldots , [9]$ of $G$ with the new label $3$.

Step 6) For every $1\leqslant j \leqslant s_k$, $9 \leqslant k \leqslant 14$ we label four edges in $\mathbf{E}_{9j}$ and $\mathbf{E}_{10j}$  with the labels $1,4$,  the  edges in $\mathbf{E}_{11j}$ and $\mathbf{E}_{12j}$  with the labels $2,4$, and the  edges in $\mathbf{E}_{13j}$ and $\mathbf{E}_{14j}$  with the labels $3,4$ as Step $1$ and $2$, respectively.

\medskip
Step 7) For every $1\leqslant j \leqslant s_k$, $15 \leqslant k \leqslant 26$ we label all four edges 
in $\mathbf{E}_{15j},\ldots ,\mathbf{E}_{18j}$ with the labels $(1,2,4,1),(1,2,4,2),(1,2,4,3),(1,2,4,4)$, the all four edges
 in $\mathbf{E}_{19j},\ldots ,\mathbf{E}_{22j}$ with the labels $(1,3,4,1),(1,3,4,2),(1,3,4,3),(1,3,4,4)$, and all the four edges
 in $\mathbf{E}_{23j},\ldots ,\mathbf{E}_{26j} $ with the labels $(2,3,4,1),(2,3,4,2),(2,3,4,3),(2,3,4,4)$, respectively.

\medskip
Step 8) For every $1\leqslant j \leqslant s_{27}$ we label the four edges in $\mathbf{E}_{27j}$ with the labels $1,2,3,4$, respectively.

So we labeled nineteen  corresponding classes  of $G$ with the new label $4$. Continuing this method we obtain that the number of  corresponding classes  of $G$ that can be labeled with the new label $m$, $m\geqslant 2$ is $2{m-1 \choose 1} +m{m-1 \choose 2}+{m-1\choose 3}$.

This labeling is distinguishing, because if $f$ is an automorphism of $G[P_2]$ preserving the labeling then there exist the automorphism $\alpha$ of $G$ and $\beta_{\alpha g}$ of $P_2$ such that  $f(g,x)=(\alpha g,\beta_{\alpha g}x)$, where $g\in V(G)$ and $x\in V(P_2)$. With respect to the method of labeling it is concluded that  $\alpha$   is the identity automorphism of $G$, because we labeled the set of four edges corresponding to the edge $e_{ij}$ of $G$, for every $1\leqslant j \leqslant s_i$ the same and different from the corresponding edges to the edge $e_{kj}$ of $G$ where $i\neq k$. On the  other hand  $\beta_{\alpha g}$ is the identity automorphism on $P_2$, because for each four edges corresponding to an edge of $G$, none of two distinct labels can not be repeated (at most, one of labels can be repeated).  Therefore $f$ is the identity automorphism of  $G[P_2]$. Since we used $min\{k : \sum_{m=2}^{k}\left( 2{m-1 \choose 1} +m{m-1 \choose 2}+{m-1\choose 3} \right)\geqslant D'(G)\}$ labels, the result follows.\qed

\begin{theorem}\label{thm G[H]}
Let $G$ and $H$ be two connected graphs with $Aut(G[H])= (Aut (G))[Aut (H)]$. If $\vert V(G) \vert \leqslant \vert E(H) \vert +1$, then $D'(G[H])\leqslant 2$.
\end{theorem}
\proof
Since $\vert V(G) \vert \leqslant \vert E(H) \vert +1$, we can label the edges  of $i$th copy of $H$ with $i-1$ labels $1$ and $ \vert E(H) \vert- (i-1)$ labels $2$, for every $1 \leqslant i \leqslant \vert V(G) \vert $. On the other hand each edge of $G$ is correspond to $\vert V(H) \vert^2 $ edges in $G[H]$. Let $V(G) =\{g_1,\ldots , g_{\vert V(G) \vert}\}$ and $V(H) =\{h_1,\ldots , h_{\vert V(H) \vert}\}$. If $e=g_ig_j$ is an edge of $G$ then $e$ is replaced by  the edges $(g_i,h_k)(g_j,h_{k'})$, where $k,k' \in \{1,\ldots , \vert V(H) \vert\}$. We label the edges $(g_i,h_p)(g_j,h_1),\ldots , (g_i,h_p)(g_j,h_{\vert V(H) \vert})$ with $p-1$ labels $2$ and $\vert V(H) \vert -(p-1)$ labels $1$, where $1 \leqslant p \leqslant \vert V(H) \vert$. We do the same work for the remaining edges of $G$.

As every copy of $H$ has a different number of edges with label $2$, they can not be interchanged. The same is true for the edges of each copy of $H$. Therefore the labeling is $2$-distinguishing labeling.\qed

\begin{corollary}\label{cor idempo}
Let $G$  be a connected graph  such that $Aut(G[G])= (Aut (G))[Aut (G)]$, then $D'(G[G])\leqslant 2$.
\end{corollary}
\proof
Since $G$ satisfies the conditions of Theorem \ref{thm G[H]}, the result follows.\qed

\begin{corollary}
Let $G$  be a connected graph  such that $Aut(G[G])=Aut(G)[Aut (G)]$ ($k\geq 2$), then $D'(G^{k})\leqslant 2$.
\end{corollary}
\proof
The proof is by induction on $k$. Let $k=1$, then the result is obtained from Corollary \ref{cor idempo}. For the induction step, we apply Theorem \ref{thmindex} by taking $H=G^{k-1}$, because $\vert V(G) \vert \leqslant \vert E(G^{k-1}) \vert +1$.\qed

\end{document}